\documentclass{birkjour}
\usepackage[utf8]{inputenc}

\usepackage{amsmath,amssymb,amsthm,blkarray,caption,esint,graphicx,hyperref,interval,mathrsfs,mathtools,multicol,orcidlink,pgfplots,tikz,tikz-cd,xcolor}

\DeclareMathOperator{\dis}{dis}
\DeclareMathOperator{\ins}{i}
\DeclareMathOperator{\ext}{e}

\DeclareMathOperator{\kernel}{ker}

\newtheorem{counter}{}[section]

\theoremstyle{definition}

\newtheorem{remark}[counter]{Remark}

\theoremstyle{plain}
\newtheorem{corollary}[counter]{Corollary}
\newtheorem{lemma}[counter]{Lemma}
\newtheorem{proposition}[counter]{Proposition}
\newtheorem{theorem}[counter]{Theorem}

\def\Xint#1{\mathchoice
   {\XXint\displaystyle\textstyle{#1}}%
   {\XXint\textstyle\scriptstyle{#1}}%
   {\XXint\scriptstyle\scriptscriptstyle{#1}}%
   {\XXint\scriptscriptstyle\scriptscriptstyle{#1}}%
   \!\int}
\def\XXint#1#2#3{{\setbox0=\hbox{$#1{#2#3}{\int}$}
     \vcenter{\hbox{$#2#3$}}\kern-.5\wd0}}

\def\dashint{\Xint-}

\renewcommand{\subjclassname}{Mathematics Subject Classification}

\begin{document}
\title{The double-layer potential for spectral constants revisited}

\author[F.L.~Schwenninger]{Felix~L.~Schwenninger\,\orcidlink{0000-0002-2030-6504}}
\address{
University of Twente\\
P.O.~Box 217\\
7500~AE Enschede\\
The Netherlands}
\email{f.l.schwenninger@utwente.nl}

\author[J.~de~Vries]{Jens~de~Vries\,\orcidlink{0009-0008-2122-814X}}
\address{
University of Twente\\
P.O.~Box 217\\
7500~AE Enschede\\
The Netherlands}
\email{j.devries-4@utwente.nl}

\begin{abstract}
We thoroughly analyse the double-layer potential's role in approaches to spectral sets in the spirit of Delyon--Delyon, Crouzeix and Crouzeix--Palencia. 
While the potential is well-studied, we aim to clarify on several of its aspects in light of these references. In particular, we illustrate how the associated integral operators can be used to characterize the convexity of the domain and the inclusion of the numerical range in its closure. We furthermore give a direct proof of a result by Putinar--Sandberg---a generalization of Berger--Stampfli's mapping theorem---circumventing dilation theory. Finally, we show for matrices that any smooth domain whose closure contains the numerical range admits a spectral constant only depending on the extremal function and vector. This constant is consistent with the so far best known absolute bound $1+\sqrt{2}$.
\end{abstract}

\thanks{The second named author is financed by the Dutch Research Council (NWO) grant OCENW.M20.292.}
\subjclass{Primary  47A25; Secondary 47A12, 47B91}
\keywords{Spectral sets, double-layer potential, numerical range.} 

\maketitle

\section*{Introduction}
Let $A$ be an operator on a Hilbert space $H$ and let $X$ be a bounded subset of $\mathbb{C}$ such that the closure $X^{-}$ contains the spectrum of $A$. A positive number $\kappa$ is called a \textit{spectral constant} for $A$ on $X$ if for every rational $f\colon X^{-}\to\mathbb{C}$ with poles off $X^{-}$ the estimate
\begin{align*}
    \|f(A)\|\leq\kappa\cdot\sup_{z\in X}|f(z)|
\end{align*}
holds. The goal is to find such spectral constants (possibly independent of $A$) for given candidates of sets $X$ (which will depend on $A$). We refer to \cite{badea2014spectral} for an extensive overview that includes many examples and references. It is  folklore that $\kappa$ is a spectral constant for $A$ on $X$ if and only if $\kappa$ is a spectral constant for $A$ on every smoothly bounded open neighbourhood $\Omega$ of the spectrum of $A$ for which $\Omega^{-}$ contains $X$. For the sake of completeness we include a proof of this fact in the appendix, see Proposition \ref{smoothApproximation}. Therefore, we may focus on smoothly bounded open sets $\Omega$ that contain the spectrum of $A$. Throughout this paper we endow $\partial\Omega$ with the orientation that is positive with respect to the outward normal vectors to $\Omega$.

Let $\mathcal{A}(\Omega^{-})$ be the space of continuous functions from $\Omega^{-}$ to $\mathbb{C}$ that are analytic on $\Omega$. The maximum modulus principle implies that $\mathcal{A}(\Omega^{-})$ embeds isometrically into $\mathcal{C}(\partial\Omega)$ with respect to the supremum norms. The \textit{analytic functional calculus} $\gamma\colon\mathcal{A}(\Omega^{-})\to\mathcal{L}(H)$ defined by the Cauchy integral
\begin{align*}
    \gamma(f):=\frac{1}{2\pi i}\int_{\partial\Omega}f(\sigma)(\sigma-A)^{-1} \ \mathrm{d}\sigma
\end{align*}
for $f\in\mathcal{A}(\Omega^{-})$ is a unital algebra homomorphism and extends the definition of the rational functional calculus. The estimate
\begin{align*}
\|\gamma(f)\|\leq\frac{1}{2\pi}\int_{\partial\Omega}\|(\sigma-A)^{-1}\| \ |\mathrm{d}\sigma|\cdot\sup_{z\in\Omega}|f(z)|
\end{align*}
for all $f\in\mathcal{A}(\Omega^{-})$ shows that $\gamma$ is bounded and yields a spectral constant that depends on both $A$ and $\Omega$. That being said, it is a priori not clear that spectral constants independent of the operator exist, but, surprisingly, they do arise  under stronger assumptions on the domain $\Omega$ than just the inclusion of the spectrum.

A key tool for finding spectral constants, initiated by Delyon--Delyon \cite{delyon1999generalization}, is the operator-valued double-layer potential, which is used to study the `symmetrized  functional calculus'
\begin{align}\label{SFC}
    f\mapsto \gamma(K_{\Omega}(f)^{*})^{*}+\gamma(f),\qquad f\in\mathcal{A}(\Omega^{-})
\end{align}
and, in particular, to bound the operator norm thereof. Here $K_{\Omega}(f)\colon\Omega^{-}\to\mathbb{C}$ is the conjugated Cauchy integral
\begin{align}\label{CCI}
		K_{\Omega}(f)(z)=\Big(\frac{1}{2\pi i}\int_{\partial\Omega}\frac{f(\sigma)^{*}}{\sigma-z} \ \mathrm{d}\sigma\Big){}^{*}
\end{align}
for all interior points $z\in\Omega$. The study of spectral constants through \eqref{SFC} originally concentrated on convex $\Omega$ \cite{delyon1999generalization,crouzeix2007numerical,crouzeix2017numerical}, but was later also expanded to non-convex $\Omega$ such as annuli \cite{caldwell2018some,crouzeix2019spectral,jury2023positivity}. 

The mapping \eqref{SFC} has been particularly influential in the investigation of spectral constants for $A$ on the numerical range $W(A)$. Crouzeix \cite{crouzeix2004bounds} conjectured that the optimal spectral constant for $A$ on $W(A)$ is at most $2$. Since $W(A)$ is convex, the smoothly bounded approximations $\Omega$ may be taken convex as well, see Proposition \ref{smoothApproximation} in the appendix. Under these assumptions---the convexity of $\Omega$ and the inclusion of $W(A)$ in $\Omega^{-}$---the operator norm of \eqref{SFC} equals $2$. Crouzeix--Palencia used this fact \cite{crouzeix2017numerical} to prove that $1+\sqrt{2}$ is a spectral constant for $A$ on $W(A)$. We also point out that, in this context, the mapping \eqref{SFC} possesses some interesting dilation theoretic properties, see \cite{putinar2005skew,hartz2021dilation}. In the recent preprint \cite{malman2024double} it was shown that the optimal spectral constant for $A$ on $W(A)$ is strictly smaller than $1+\sqrt{2}$ by an indefinite amount.

On a function theoretic level, the subtleties of the conjecture also become apparent for example when restricting to the case where the operators are compressions of the shift with finite Blaschke products, \cite{bickel2020crouzeix,bickel2018compressions}, also see \cite{bickel2020numerical,bickel2023crouzeix,gorkin2024three}.

In \cite{ransford2018remarks} the main ingredients of the Delyon--Delyon/Crouzeix/Crouzeix--Palencia approach were classified into `hard facts' of the double-layer potential and some abstract operator theory for norm bounds of homomorphisms. In our  recent contribution \cite{schwenninger2024abstract}, we tried to unify the latter part, leading to a streamlined reasoning for several earlier results \cite{neumann1950spektraltheorie,okubo1975constants,crouzeix2017numerical,crouzeix2019spectral}. The goal of the current contribution is to precisely analyse the ingredients of the former part---the `hard facts' of the double-layer potential---and thereby to complement \cite{schwenninger2024abstract}. By doing so we hope to reveal more of the underlying structure of {arguments involving the double-layer potential} and to clarify {questions like whether these arguments can expected to be sharpened.} 
In order to achieve this, we have to carefully examine the double-layer potential and the associated integral operators. In this paper we focus on smooth domains, but we are aware of the fact that this assumption is not necessary, see e.g.\ {\cite{crouzeix2019spectral,malman2024double}}. The advantage of assuming smoothness is that many properties of the double-layer potential quickly follow from Plemelj--Sokhotski's formulae for Cauchy transforms. Of course, we are not the first to stress the importance of the double-layer potential in the context of spectral constants, \cite{delyon1999generalization, putinar2005skew,crouzeix2007numerical,crouzeix2017numerical}. Yet, there is considerable interest in revisiting the double-layer potential as it could lead to deeper insights. In particular, we outline how the double-layer transform can be used to characterize the convexity of $\Omega$ as well as the inclusion of $W(A)$ in $\Omega^{-}$. Note that the double-layer potential naturally relates to operators  defined on the continuous functions, whereas spectral constants relate to the norm of the homomorphism acting on holomorphic.  More specifically, it is the interplay of the \emph{interior} and the \emph{singular} \emph{double-layer transform}, see below, which encodes these non-trivial discrepancy. See also \cite{malman2024double}, where this subtlety is crucial.

Many results concerning bounds for analytic functional calculi, especially on the disk, have historically been obtained through dilation theory. A celebrated example is the mapping theorem by Berger--Stapfli \cite{berger1967mapping} for numerical ranges. Using dilation theory they proved that $W(A)\subseteq\mathbb{D}^{-}$ implies $W(\gamma(f))\subseteq\mathbb{D}^{-}$ provided that $\|f\|_{\infty}\leq1$ and $f(0)=0$. Note that, independently of Berger--Stampfli, Kato \cite{kato1965some} proved the same result using a more direct method, without dilation theory. Another such proof is due to Klaja--Mashreghi--Ransford \cite{klaja2016mapping} and is based on finite Blaschke products. {Dilation techniques, similar to that of Berger--Stampfli, were later employed by Putinar--Sandberg \cite{putinar2005skew} to obtain a more general result.} They proved (for convex $\Omega$) that $W(A)\subseteq\Omega^{-}$ implies $W(\gamma(f))\subseteq\mathbb{D}^{-}$ whenever $\|f\|_{\infty}\leq1$ and \eqref{CCI} satisfies $K_{\Omega}(f)=0$. In the current paper we exploit properties of the double-layer potential to provide a self-contained proof of the result by Putinar--Sandberg, without dilation theory.

The paper is organized so that Sections \ref{sectionCauchy} and \ref{sectionDouble} focus on the function theory related to the double-layer potential $P_{\Omega}\colon(\partial\Omega\times\mathbb{C})\setminus\Delta\to\mathbb{R}$, while Section \ref{sectionOperator} discusses the operator-valued counterparts. In a broad sense, Sections \ref{sectionPutinar}, \ref{sectionRole} and \ref{sectionSpectral} explore applications of the theory in the framework of spectral constants.

In Section \ref{sectionCauchy} we consider H\"older continuous $\varphi\colon\partial\Omega\to\mathbb{C}$ and look at the restrictions of the Cauchy transform
\begin{align*}
	z\mapsto\frac{1}{2\pi i}\int_{\partial\Omega}\frac{\varphi(\sigma)}{\sigma-z} \ \mathrm{d}\sigma,\qquad z\in\mathbb{C}\setminus\partial\Omega
\end{align*}
to the interior $\Omega$ and the exterior $\mathbb{C}\setminus\Omega^{-}$. We discuss their continuous extensions to the respective closures $\Omega^{-}$ and $\mathbb{C}\setminus\Omega$, as well as their mutual `jump relations' on the shared boundary $\partial\Omega$.

In Section \ref{sectionDouble} we analyze how the previous theory can be used to describe the boundary behaviour of the restrictions of the double-layer transform
\begin{align*}
	z\mapsto\int_{\partial\Omega}\varphi(\sigma)P_{\Omega}(\sigma,z) \ |\mathrm{d}\sigma|,\qquad z\in\mathbb{C}\setminus\partial\Omega
\end{align*}
to $\Omega$ and $\mathbb{C}\setminus\Omega^{-}$. Specifically, we discuss how the Neumann--Poincar\'e operator $K_{\Omega}\colon\mathcal{C}(\partial\Omega)\to\mathcal{C}(\partial\Omega)$---a compact operator that historically has been instrumental in establishing solvability of Laplace's equation with continuous boundary data---forms a link between boundary values of both restrictions. Among other characterizations we review how the value of $\|K_{\Omega}\|$ characterizes the convexity of $\Omega$, see also \cite{kral1980integral}.

In Section \ref{sectionOperator} we consider the operator-valued double-layer potential. We relate it to the harmonic functional calculus and recover the mapping \eqref{SFC}. Under the assumption that $\Omega$ is convex we discuss various characterizations of the inclusion of $W(A)$ in $\Omega^{-}$.

In Section \ref{sectionPutinar} we illustrate how properties of \eqref{SFC} can be used to prove the theorem by Putinar--Sandberg. We present a proof that circumvents the dilation theoretic arguments.

In Section \ref{sectionRole} we explain how the operator-valued double-layer potential was used in the works \cite{delyon1999generalization,crouzeix2007numerical,crouzeix2017numerical,ransford2018remarks} to find spectral constants for (smooth approximations of) the numerical range. For each of these results we describe the explicit role of \eqref{SFC}.

In Section \ref{sectionSpectral} we use techniques from \cite{schwenninger2024abstract} to prove that the inclusion $W(A)\subseteq\Omega^{-}$ implies the bound
\begin{align*}
\|\gamma\|\leq 1+\sqrt{1-\rho(f_{0},x_{0})}
\end{align*}
for some constant $\rho(f_{0},x_{0})\in\interval{-1}{1}$ that depends on an extremal pair $f_{0}\in\mathcal{A}(\Omega^{-})$ and $x_{0}\in H$.

Clearly, topics centered around Crouzeix's conjecture span more theory than the approach discussed above and in the rest of the paper. We chose to mention those references that are closely related to the Delyon--Delyon/Crouzeix/Crouzeix--Palencia approach, which particularly focuses on general matrices or operators. This direction based on Cauchy integrals is in fact also relevant in the context of applications in numerical analysis, such as error bounds for Krylov type methods \cite{crouzeix2019spectral} and stability of Runge--Kutta schemes \cite{tadmor2023runge}. Nevertheless, there are of course interesting developments for specific classes of matrices based on completely different techniques. We refer e.g.\ to \cite{bickel2020numerical,gorkin2024three} and the references therein as well as to the recent paper \cite{oloughlin2024crouzeix}.

\section{The Cauchy transforms}\label{sectionCauchy}
Let $\Omega$ be a smoothly bounded open subset of $\mathbb{C}$. A function $\varphi\colon\partial\Omega\to\mathbb{C}$ is \textit{H\"older continuous} if there exists constants $c>0$ and $\alpha>0$ such that
\begin{align*}
	|\varphi(\sigma)-\varphi(\tau)|\leq c|\sigma-\tau|^{\alpha}
\end{align*}
for all $\sigma,\tau\in\partial\Omega$. The space of all such H\"older continuous functions, which we denote by $\mathcal{C}_{\bullet}(\partial\Omega)$, is a uniformly dense subalgebra of $\mathcal{C}(\partial\Omega)$.

Now suppose that $\varphi\in\mathcal{C}_{\bullet}(\partial\Omega)$ is given. The assignment
\begin{align*}
	z\mapsto\frac{1}{2\pi i}\int_{\partial\Omega}\frac{\varphi(\sigma)}{\sigma-z} \ \mathrm{d}\sigma,\qquad z\in\mathbb{C}\setminus\partial\Omega
\end{align*}
clearly determines two analytic functions $\Phi_{\Omega}^{\ins}(\varphi)\colon\Omega\to\mathbb{C}$ and $\Phi_{\Omega}^{\ext}(\varphi)\colon\mathbb{C}\setminus\Omega^{-}\to\mathbb{C}$. Employing principle-value integration allows one to define a similar function on the boundary. Indeed, under the H\"older continuity assumption it is known that
\begin{align*}
	z\mapsto\frac{1}{2\pi i}\dashint_{\partial\Omega}\frac{\varphi(\sigma)}{\sigma-z} \ \mathrm{d}\sigma,\qquad z\in\partial\Omega
\end{align*}
establishes a continuous function $\Phi_{\Omega}(\varphi)\colon\partial\Omega\to\mathbb{C}$, see e.g.\ \cite[Section 2.12]{muskhelishvili2008singular}. The functions $\Phi_{\Omega}^{\ins}(\varphi)$, $\Phi_{\Omega}^{\ext}(\varphi)$ and $\Phi_{\Omega}(\varphi)$ are related through \textit{Plemelj--Sokhotski's formulae}, which state that $\Phi_{\Omega}^{\ins}(\varphi)$ and $\Phi_{\Omega}^{\ext}(\varphi)$ extend continuously to $\Omega^{-}$ and $\mathbb{C}\setminus\Omega$, respectively, with boundary values
\begin{align*}
	\Phi_{\Omega}^{\ins}(\varphi)(z)=\Phi_{\Omega}(\varphi)(z)+\frac{\varphi(z)}{2},\qquad\Phi_{\Omega}^{\ext}(\varphi)(z)=\Phi_{\Omega}(\varphi)(z)-\frac{\varphi(z)}{2}
\end{align*}
for all $z\in\partial\Omega$, see e.g.\ \cite[Section 2.17]{muskhelishvili2008singular}. In particular, we trivially obtain the \textit{jump formula}
\begin{align*}
	\Phi_{\Omega}^{\ins}(\varphi)(z)-\Phi_{\Omega}^{\ext}(\varphi)(z)=\varphi(z)
\end{align*}
for all $z\in\partial\Omega$. 

\begin{remark}
		The H\"older continuity is necessary. Let $\varphi\colon\partial\mathbb{D}\to\mathbb{C}$ be a (not necessarily H\"older) continuous function with Fourier coefficients $(c_{k})_{k\in\mathbb{Z}}$. Define $\Phi_{\mathbb{D}}^{\ins}(\varphi)\colon\mathbb{D}\to\mathbb{C}$ as above and let $r\in\interval[open right]{0}{1}$ be arbitrary. An application of the dominated convergence theorem yields
		\begin{align*}
			\Phi_{\mathbb{D}}^{\ins}(\varphi)(rz)&=\frac{1}{2\pi i}\int_{\partial\mathbb{D}}\frac{\varphi(\sigma)}{\sigma-rz} \ \mathrm{d}\sigma=\frac{1}{2\pi i}\int_{\partial\mathbb{D}}\frac{\varphi(\sigma)}{\sigma}\sum_{k=0}^{\infty}\frac{r^{k}z^{k}}{\sigma^{k}} \ \mathrm{d}\sigma\\
			&=\sum_{k=0}^{\infty}\frac{1}{2\pi }\int_{\partial\mathbb{D}}\frac{\varphi(\sigma)}{\sigma^{k}} \ |\mathrm{d}\sigma| \ r^{k}z^{k}=\sum_{k=0}^{\infty}c_{k}r^{k}z^{k}
		\end{align*}
		for all $z\in\partial\mathbb{D}$. Let $(e_{k})_{k\in\mathbb{Z}}$ be the Fourier basis of $L^{2}(\partial\mathbb{D})$ and write $\varphi_{+}$ for the image of $\varphi$ under the orthogonal projector of $L^{2}(\partial\mathbb{D})$ onto the linear span of $\{e_{0},e_{1},\ldots\}$. In other words, $\varphi_{+}$ is the Riesz projection of $\varphi$. By elementary Hardy space theory the radial limits satisfy
		\begin{align*}
			\lim_{r\nearrow1}\Phi_{\mathbb{D}}^{\ins}(\varphi)(rz)=\lim_{r\nearrow1}\sum_{k=0}^{\infty}c_{k}r^{k}z^{k}=\varphi_{+}(z)
		\end{align*}
		for almost all $z\in\partial\mathbb{D}$. In particular, a continuous extension of $\Phi_{\mathbb{D}}^{\ins}(\varphi)$ to $\mathbb{D}^{-}$, if it exists, must coincide with $\varphi_{+}$ on $\partial\mathbb{D}$. Since there are continuous functions $\varphi$ for which $\varphi_{+}$ is not continuous, such a continuous extension need not exist in general, see e.g.\ \cite[Example 6.2]{pohl2009advanced}.
	\end{remark}

For a closed subset $C$ of $\mathbb{C}$ with non-empty boundary $\partial C$ we write $\mathcal{A}(C)$ for the space of continuous functions from $C$ to $\mathbb{C}$ that are analytic on the interior of $C$ and vanish at infinity if $C$ is unbounded. The maximum modulus principle implies that $\mathcal{A}(C)$ embeds isometrically into $\mathcal{C}(\partial C)$ with respect to the supremum norms. We obtain three well-defined operators.
\begin{itemize}
    \item The \textit{interior Cauchy transform} $\Phi_{\Omega}^{\ins}\colon\mathcal{C}_{\bullet}(\partial\Omega)\to\mathcal{A}(\Omega^{-})$.
    \item The \textit{exterior Cauchy transform} $\Phi_{\Omega}^{\ext}\colon\mathcal{C}_{\bullet}(\partial\Omega)\to\mathcal{A}(\mathbb{C}\setminus\Omega)$.
    \item The \textit{singular Cauchy transform} $\Phi_{\Omega}\colon\mathcal{C}_{\bullet}(\partial\Omega)\to\mathcal{C}(\partial\Omega)$.
\end{itemize}
The values of the Cauchy transforms are in fact again H\"older continuous on the boundary, see e.g.\ \cite[Sections 2.19 and 2.20]{muskhelishvili2008singular}.

\section{The double-layer transforms}\label{sectionDouble}
Let $\Omega$ be a smoothly bounded open subset of $\mathbb{C}$. Let $n_{\Omega}\colon\partial\Omega\to\mathbb{C}$ be the function that sends any boundary point in $\partial\Omega$ to the outward unit normal vector of $\Omega$ at that point. Suppose that $\Delta$ is the diagonal in $\partial\Omega\times\partial\Omega$ and define the \textit{double-layer potential} $P_{\Omega}\colon(\partial\Omega\times\mathbb{C})\setminus\Delta\to\mathbb{R}$ by
\begin{align*}
	P_{\Omega}(\sigma,z):=\frac{1}{\pi}\Re\Big(\frac{n_{\Omega}(\sigma)}{\sigma-z}\Big)
\end{align*}
for $\sigma\in\partial\Omega$ and $z\in\mathbb{C}$ with $\sigma\neq z$.

Let $\varphi\in\mathcal{C}_{\bullet}(\Omega^{-})$ be given. The assignment
\begin{align*}
	z\mapsto\int_{\partial\Omega}\varphi(\sigma)P_{\Omega}(\sigma,z) \ |\mathrm{d}\sigma|,\qquad z\in\mathbb{C}\setminus\partial\Omega
\end{align*}
yields two harmonic functions $K_{\Omega}^{\ins}(\varphi)\colon\Omega\to\mathbb{C}$ and $K_{\Omega}^{\ext}(\varphi)\colon\mathbb{C}\setminus\Omega^{-}\to\mathbb{C}$ as the former is the sum of $\Phi_{\Omega}^{\ins}(\varphi)$ and $\Phi_{\Omega}^{\ins}(\varphi^{*})^{*}$ on $\Omega$ and the latter is the sum of $\Phi_{\Omega}^{\ext}(\varphi)$ and $\Phi_{\Omega}^{\ext}(\varphi^{*})^{*}$ on $\mathbb{C}\setminus\Omega^{-}$. Similarly, the boundary version
\begin{align*}
	z\mapsto\dashint_{\partial\Omega}\varphi(\sigma)P_{\Omega}(\sigma,z) \ |\mathrm{d}\sigma|,\qquad z\in\partial\Omega
\end{align*}
yields a continuous function $K_{\Omega}(\varphi)\colon\partial\Omega\to\mathbb{C}$ as it is the sum of $\Phi_{\Omega}(\varphi)$ and $\Phi_{\Omega}(\varphi^{*})^{*}$ on $\partial\Omega$. A straightforward application of Plemelj--Sokhotski's formulae shows that $K_{\Omega}^{\ins}(\varphi)$ and $K_{\Omega}^{\ext}(\varphi)$ extend continuously to $\Omega^{-}$ and $\mathbb{C}\setminus\Omega$, respectively, with boundary values
\begin{align*}
	K_{\Omega}^{\ins}(\varphi)(z)=K_{\Omega}(\varphi)(z)+\varphi(z),\qquad K_{\Omega}^{\ext}(\varphi)(z)=K_{\Omega}(\varphi)(z)-\varphi(z)
\end{align*}
for all $z\in\partial\Omega$.

For a closed subset $C$ of $\mathbb{C}$ with non-empty boundary $\partial C$ we write $\mathcal{H}(C)$ for the uniform closure of $\mathcal{A}(C)+\mathcal{A}(C)^{*}$ in $\mathcal{C}(\partial C)$. We obtain three well-defined operators.
\begin{itemize}
    \item The \textit{interior double-layer transform} $K_{\Omega}^{\ins}\colon\mathcal{C}_{\bullet}(\partial\Omega)\to\mathcal{H}(\Omega^{-})$.
    \item The \textit{exterior double-layer transform} $K_{\Omega}^{\ext}\colon\mathcal{C}_{\bullet}(\partial\Omega)\to\mathcal{H}(\mathbb{C}\setminus\Omega)$.
    \item The \textit{singular double-layer transform} $K_{\Omega}\colon\mathcal{C}_{\bullet}(\partial\Omega)\to\mathcal{C}(\partial\Omega)$.
\end{itemize}
The singular double-layer transform is also known as the \textit{Neumann--Poincar\'e operator}. In contrast to the Cauchy transforms, the double-layer transforms are bounded. To see this we need following technical lemma, see e.g.\ \cite[Lemma 3.15]{folland2020introduction}. 
\begin{lemma}\label{boundOperatorNormK}
    There exists a $\delta>0$ such that 
	\begin{align*}
		|\Re((\sigma-\tau)^{*}n_{\Omega}(\sigma))|\leq\delta|\sigma-\tau|^{2}
	\end{align*}
	for all boundary points $\sigma,\tau\in\partial\Omega$.
\end{lemma}
\begin{proof}
	It suffices to prove the estimate on some neighbourhood of the diagonal $\Delta$ in $\partial\Omega\times\partial\Omega$. By compactness of $\partial\Omega$ there are finitely many open intervals $I_{1},\ldots,I_{l}$ in $\mathbb{R}$ and open sets $D_{1},\ldots,D_{l}$ in $\mathbb{C}$ such that
	\begin{align*}
		\Delta\subseteq\bigcup_{k=1,\ldots,l}(\partial\Omega\cap D_{k})\times(\partial\Omega\cap D_{k}).
	\end{align*}
	and for each $k=1,\ldots,l$ the set $\partial\Omega\cap D_{k}$ equals the graph of some smooth function $F_{k}\colon I_{k}\to\mathbb{R}$ with bounded second derivative.	Without loss of generality we may assume that
	\begin{align*}
		\partial\Omega\cap D_{k}=\{\xi+iF_{k}(\xi):\xi\in I_{k}\}.
	\end{align*}
	Suppose that $\tfrac{1}{2}|F_{k}''(\xi)|\leq\delta_{k}$ for all $\xi\in I_{k}$. If $\sigma=s+iF_{k}(s)$ and $\tau=t+iF_{k}(t)$ for some $s,t\in I_{k}$, then it follows from direct computation that
	\begin{align*}
		|\Re((\sigma-\tau)^{*}n_{\Omega}(\sigma))|=\Big|\frac{(s-t)F_{k}'(s)-(F_{k}(s)-F_{k}(t))}{\sqrt{1+F_{k}'(s)^{2}}}\Big|\leq\delta_{k}|s-t|^{2}\leq\delta_{k}|\sigma-\tau|^{2}
	\end{align*}
	by the Taylor expansion theorem. Now take $\delta:=\max\{\delta_{1},\ldots,\delta_{k}\}$.
\end{proof}
It follows from Lemma \ref{boundOperatorNormK} that 
\begin{align*}
K_{\Omega}(\varphi)(z)=\int_{\partial\Omega}\varphi(\sigma)P_{\Omega}(\sigma,z) \ |\mathrm{d}\sigma|
\end{align*}
for all $\varphi\in\mathcal{C}_{\bullet}(\partial\Omega)$ and $z\in\partial\Omega$ as an ordinary integral, that is, principle-value integration is not necessary to deal with the singularity of the integrand at $\sigma=z$. Moreover, we find a $\delta>0$ such that 
\begin{align*}
	\|K_{\Omega}(\varphi)\|_{\infty}\leq\frac{\delta}{\pi}|\partial\Omega|\|\varphi\|_{\infty}
\end{align*}
for all $\varphi\in\mathcal{C}_{\bullet}(\partial\Omega)$. Thus $K_{\Omega}^{\ins}$, $K_{\Omega}^{\ext}$ and $K_{\Omega}$ are bounded and admit bounded extensions to $\mathcal{C}(\partial\Omega)$. 
\begin{remark}
    It follows from \cite[Theorem 2.30]{kress2013linear} that $K_{\Omega}$ is a compact self-mapping on $\mathcal{C}(\partial\Omega)$.
\end{remark}

\subsection{Convexity of $\Omega$}
The convexity of $\Omega$ is encoded by a positivity condition on the double-layer potential $P_{\Omega}$. The following result seems to be folklore.
\begin{proposition}\label{characterizationConvexity1} 
 Let $\Delta$ be the diagonal in $\partial\Omega\times\partial\Omega$. The domain $\Omega$ is convex if and only if the function $P_{\Omega}$ is positive on $(\partial\Omega\times\partial\Omega)\setminus\Delta$.
\end{proposition}
\begin{proof}
	For each point $\sigma$ on the boundary $\partial\Omega$ let $\Pi_{\sigma}$ be the unique closed half-plane in $\mathbb{C}$ with $n_{\Omega}(\sigma)$ as outward normal vector at $\sigma$. The domain $\Omega$ is convex if and only if
	\begin{align*}    \partial\Omega\subseteq\bigcap_{\sigma\in\partial\Omega}\Pi_{\sigma}.
	\end{align*}
For any $z\in\mathbb{C}$ with $\sigma\neq z$ the equality
\begin{align*}
		\Re((\sigma-z)^{*}n_{\Omega}(\sigma))=\pi|\sigma-z|^{2}P_{\Omega}(\sigma,z)
	\end{align*}
holds, implying that $z\in\Pi_{\sigma}$ if and only if $P_{\Omega}(\sigma,z)\geq0$. The result follows.
\end{proof}
The property whether $\Omega$ is convex is also encoded by the operator norm of $K_{\Omega}$. This result, for which we present a short proof, traces back to \cite[Theorem 3.1]{kral1980integral}.

\begin{proposition}\label{characterizationConvexity2}
        The following assertions are equivalent:
	\begin{itemize}
		\item[{\normalfont(i)}] $\Omega$ is a convex set,
		\item[{\normalfont(ii)}] $K_{\Omega}$ is a positive operator,
		\item[{\normalfont(iii)}] the operator norm of $K_{\Omega}$ equals $1$.
	\end{itemize}
\end{proposition}
\begin{proof}
	Let $\Delta$ be the diagonal in $\partial\Omega\times\partial\Omega$. Note that, by Proposition \ref{characterizationConvexity1}, $\Omega$ is convex if and only if $P_{\Omega}$ is positive on $(\partial\Omega\times\partial\Omega)\setminus\Delta$. The latter is true if and only if $K_{\Omega}$ is positive. Since $K_{\Omega}$ is a unital operator between unital C*-algebras, it is well-known that $K_{\Omega}$ is positive if and only if $K_{\Omega}$ is contractive, see e.g.\ \cite[Chapter 2]{paulsen2002completely}.
\end{proof}
\begin{remark}
    If $\Omega$ is convex, then the operator norm of $K_{\Omega}|\mathcal{A}(\Omega^{-})$ equals $1$ as well. In \cite[Section 5]{crouzeix2019spectral} the annulus is given as a counter-example to the converse of this statement.
\end{remark}
\subsection{The restriction of $K_{\Omega}$ to $\mathcal{A}(\Omega^{-})$}
Let $\Omega$ be a smoothly bounded open subset of $\mathbb{C}$. In this section we consider $K_{\Omega}|\mathcal{A}(\Omega^{-})$ and discuss some of its properties. First of all, it is closely related to the interior Cauchy transform. 
\begin{proposition}\label{conjugatedInteriorCauchyTransform}
The singular double-layer potential $K_{\Omega}$ maps $\mathcal{A}(\Omega^{-})$ into $\mathcal{A}(\Omega^{-})^{*}$ with values
\begin{align*}
	K_{\Omega}(f)=\Phi_{\Omega}^{\ins}(f^{*})^{*}
\end{align*}
for all $f\in\mathcal{A}(\Omega^{-})\cap\mathcal{C}_{\bullet}(\partial\Omega)$. 
\end{proposition}
\begin{proof}
For $f\in\mathcal{A}(\Omega^{-})\cap\mathcal{C}_{\bullet}(\partial\Omega)$ it holds that $\Phi_{\Omega}^{\ins}(f)=f$ and therefore
\begin{align*}
	K_{\Omega}(f)=K_{\Omega}^{\ins}(f)-f=\Phi_{\Omega}^{\ins}(f)+\Phi_{\Omega}^{\ins}(f^{*})^{*}-f=\Phi_{\Omega}^{\ins}(f^{*})^{*}\in\mathcal{A}(\Omega^{-})^{*}.
\end{align*}
The result follows as $\mathcal{A}(\Omega^{-})\cap\mathcal{C}_{\bullet}(\partial\Omega)$ is uniformly dense in $\mathcal{A}(\Omega^{-})$.
\end{proof}
In particular, Proposition \ref{conjugatedInteriorCauchyTransform} yields $K_{\Omega}(1)=1$ and consequently $K_{\Omega}^{\ins}(1)=2$ and $K_{\Omega}^{\ext}(1)=0$. Equivalently,
\begin{align*}
\int_{\partial\Omega}P_{\Omega}(\sigma,z) \ |\mathrm{d}\sigma|=\begin{cases}
    1&z\in\partial\Omega\\
    2&z\in\Omega\\
    0&z\in\mathbb{C}\setminus\Omega^{-}
\end{cases}
\end{align*}
for all $z\in\mathbb{C}$. Moreover, by approximation we obtain
	\begin{align*}
		K_{\Omega}(f)(z)=\Big(\frac{1}{2\pi i}\int_{\partial\Omega}\frac{f(\sigma)^{*}}{\sigma-z} \ \mathrm{d}\sigma\Big){}^{*}
	\end{align*}
	for all $f\in\mathcal{A}(\Omega^{-})$ and interior points $z\in\Omega$.

The following result describes the kernel of $K_{\Omega}|\mathcal{A}(\Omega^{-})$. 
\begin{proposition}\label{kernelK}
 It holds that
	\begin{align*}
		\kernel(K_{\Omega}|\mathcal{A}(\Omega^{-}))=\mathcal{A}(\Omega^{-})\cap\mathcal{A}(\mathbb{C}\setminus\Omega)^{*}.
	\end{align*}
	In particular, $\kernel(K_{\Omega}|\mathcal{A}(\Omega^{-}))$ is a subalgebra of $\mathcal{A}(\Omega^{-})$.
\end{proposition}
\begin{proof}
	Let $f\in\mathcal{A}(\Omega^{-})$ be arbitrary. On the one hand, if $K_{\Omega}(f)=0$ and $(f_{n})_{n\in\mathbb{N}}$ is a sequence in $\mathcal{A}(\Omega^{-})\cap \mathcal{C}_{\bullet}(\partial\Omega)$ that converges uniformly to $f$, then $(\Phi_{\Omega}^{\ins}(f_{n}^{*})^{*})_{n\in\mathbb{N}}$ converges uniformly to $0$ by Proposition \ref{conjugatedInteriorCauchyTransform}	and therefore
	\begin{align*}
		f=\lim_{n\to\infty}\Phi_{\Omega}^{\ins}(f_{n}^{*})^{*}-\Phi_{\Omega}^{\ext}(f_{n}^{*})^{*}=\lim_{n\to\infty}-\Phi_{\Omega}^{\ext}(f_{n}^{*})^{*}\in\mathcal{A}(\mathbb{C}\setminus\Omega)^{*}
	\end{align*}
	by the jump formula. On the other hand, if $f\in\mathcal{A}(\mathbb{C}\setminus\Omega)^{*}$, then clearly $K_{\Omega}(f)=0$ by Cauchy's integral formula. 
\end{proof}
In other words, Proposition \ref{kernelK} says that $\ker(K_{\Omega}|\mathcal{A}(\Omega^{-}))$ consists precisely of those functions in $\mathcal{A}(\Omega^{-})$ that admit an anti-analytic extension to the exterior $\mathbb{C}\setminus\Omega^{-}$ and vanish at infinity. The question whether $\mathcal{A}(\Omega^{-})$ and $\mathcal{A}(\mathbb{C}\setminus\Omega)^{*}$ intersect non-trivially is often referred to as the \textit{matching problem}. The following elementary result solves the matching problem for open disks, see \cite[Subsection 7.2.5]{remmert2012theory}. For convenience of the reader we include the simple proof.
\begin{proposition}\label{valuesDiskCaseK}
	If $\Omega$ is an open disk in $\mathbb{C}$ with center $c$, then
	\begin{align*}
		K_{\Omega}(f)=f(c)
	\end{align*}
	for all $f\in\mathcal{A}(\Omega^{-})$.
\end{proposition}
\begin{proof}
	Suppose that $|\sigma-c|=r$ for all $\sigma\in\partial\Omega$. For each $z\in\Omega$ the expression
	\begin{align*}
		K_{\Omega}(f)(z)-f(c)&=\Big(\frac{1}{2\pi i}\int_{\partial\Omega}\frac{f(\sigma)^{*}}{\sigma-z} \ \mathrm{d}\sigma\Big){}^{*}-\frac{1}{2\pi i}\int_{\partial\Omega}\frac{f(\sigma)}{\sigma-c} \ \mathrm{d}\sigma\\
		&=\Big(\frac{1}{2\pi r}\int_{\partial\Omega}\frac{f(\sigma)^{*}(\sigma-c)}{\sigma-z} \ |\mathrm{d}\sigma|\Big){}^{*}-\frac{1}{2\pi r}\int_{\partial\Omega}f(\sigma) \ |\mathrm{d}\sigma|\\
		&=\frac{1}{2\pi r}\int_{\partial\Omega}\frac{f(\sigma)(z-c)^{*}}{(\sigma-c)^{*}-(z-c)^{*}} \ |\mathrm{d}\sigma|\\
		&=\frac{1}{2\pi i}\int_{\partial\Omega}\frac{f(\sigma)(z-c)^{*}}{r^{2}-(\sigma-c)(z-c)^{*}} \ \mathrm{d}\sigma
	\end{align*}
	vanishes by Cauchy's integral theorem as the continuous integrand
	\begin{align*}
		\sigma\mapsto\frac{f(\sigma)(z-c)^{*}}{r^{2}-(\sigma-c)(z-c)^{*}},\qquad\sigma\in\Omega^{-}
	\end{align*}
	is analytic on $\Omega$.
\end{proof}
Indeed, Proposition \ref{kernelK} and Proposition \ref{valuesDiskCaseK} imply that, when $\Omega$ is an open disk, $\mathcal{A}(\Omega^{-})\cap\mathcal{A}(\mathbb{C}\setminus\Omega)^{*}$ consists of the functions in $\mathcal{A}(\Omega^{-})$ that vanish at the center of $\Omega$. Interestingly, the matching problem for a domain whose boundary is a non-circular ellipse has a negative answer, see \cite[Theorem 3.4]{ebenfelt2001inverse}.
\section{The operator-valued double-layer potential}\label{sectionOperator}
Let $A$ be a bounded operator on a Hilbert space $H$. Let $\Omega$ be a smoothly bounded open subset of $\mathbb{C}$ that contains the spectrum of $A$. For $\sigma\in\partial\Omega$ we consider the Hermitian operator
\begin{align*}
	P_{\Omega}(\sigma,A):=\frac{1}{\pi}\Re(n_{\Omega}(\sigma)(\sigma-A)^{-1}).
\end{align*}

\subsection{The harmonic functional calculus}
Moving forward, we turn our attention to $\mathcal{H}(\Omega^{-})$.
\begin{remark}
The spaces $\mathcal{H}(\Omega^{-})$ and $\mathcal{C}(\partial\Omega)$ are equal if and only if $\mathbb{C}\setminus\Omega^{-}$ is connected, see e.g.\ \cite[Section 6.5]{conway1991theory}.
\end{remark}
In the remainder of this section we assume that $\Omega$ is connected. In this case the analytic functional calculus $\gamma \colon\mathcal{A}(\Omega^{-})\to\mathcal{L}(H)$ of $A$ on $\Omega$ can be algebraically extended to the harmonic setting. Suppose that $f_{1},g_{1},f_{2},g_{2}\in\mathcal{A}(\Omega^{-})$ satisfy $f_{1}+g_{1}^{*}=f_{2}+g_{2}^{*}$. Connectedness implies that $f_{1}-f_{2}=\lambda$ and $g_{1}-g_{2}=-\lambda^{*}$ for some $\lambda\in\mathbb{C}$ and therefore
\begin{align*}
(\gamma (f_{1})+\gamma (g_{1})^{*})-(\gamma (f_{2})+\gamma (g_{2})^{*})=\gamma (f_{1}-f_{2})+\gamma (g_{1}-g_{2})^{*}=\lambda+(-\lambda^{*})^{*}=0.
\end{align*}
Consequently, we obtain a well-defined bounded operator $\hat{\gamma}\colon\mathcal{H}(\Omega^{-})\to\mathcal{L}(H)$, called the \textit{harmonic functional calculus} of $A$ on $\Omega$, that is densely defined by
\begin{align*}
	\hat{\gamma}(f+g^{*}):=\gamma (f)+\gamma (g)^{*}
\end{align*}
for $f,g\in\mathcal{A}(\Omega^{-})$. It is clear that $\hat{\gamma}$ extends $\gamma$. Next we consider the composition
\begin{center}
\begin{tikzcd}
    \mathcal{C}(\partial\Omega)\arrow[rr,"K_{\Omega}^{\ins}"]&&\mathcal{H}(\Omega^{-})\arrow[rr,"\hat{\gamma}"]&&\mathcal{L}(H).
\end{tikzcd}
\end{center}
Let us first break down the action of $\hat{\gamma}\circ K_{\Omega}^{\ins}$ on $\mathcal{A}(\Omega^{-})$.
\begin{proposition}\label{decomposition}
    Assume that $\Omega$ is connected. It holds that
    \begin{align*}
        (\hat{\gamma}\circ K_{\Omega}^{\ins})(f)=\gamma(K_{\Omega}(f)^{*})^{*}+\gamma(f)
    \end{align*}
    for all $f\in\mathcal{A}(\Omega^{-})$.
\end{proposition}
\begin{proof}
Recall that $K_{\Omega}^{\ins}(\varphi)=K_{\Omega}(\varphi)+\varphi$ for all $\varphi\in\mathcal{C}(\partial\Omega)$ by Plemelj--Sokhotski's formulae. By Proposition \ref{conjugatedInteriorCauchyTransform} we have $K_{\Omega}(f)\in\mathcal{A}(\Omega^{-})^{*}$ for all $f\in\mathcal{A}(\Omega^{-})$. The result follows by definition of the harmonic functional calculus.
\end{proof}
Next we observe that $\hat{\gamma}\circ K_{\Omega}^{\ins}$ has an integral representation in terms of the operator-valued double-layer potential.
\begin{proposition}\label{integralRepresentationOperatorValuedDP}
	Assume that $\Omega$ is connected. The equality
	\begin{align*}
		(\hat{\gamma}\circ K_{\Omega}^{\ins})(\varphi)=\int_{\partial\Omega}\varphi(\sigma)P_{\Omega}(\sigma,A) \ |\mathrm{d}\sigma|
	\end{align*}
	holds for all $\varphi\in\mathcal{C}(\partial\Omega)$.
\end{proposition}
\begin{proof}
	It suffices to prove the equality for $\varphi\in\mathcal{C}_{\bullet}(\partial\Omega)$. Let $D$ be a smoothly bounded open subset of $\mathbb{C}$ such that $\Omega$ contains $D^{-}$ and $D$ contains the spectrum of $A$. We have
	\begin{align*}
		\gamma (\Phi_{\Omega}^{\ins}(\varphi))&=\frac{1}{2\pi i}\int_{\partial D}\Phi_{\Omega}^{\ins}(\varphi)(\tau)(\tau-A)^{-1} \ \mathrm{d}\tau\\
		&=\frac{1}{2\pi i}\int_{\partial D}\Big(\frac{1}{2\pi i}\int_{\partial\Omega}\frac{\varphi(\sigma)}{\sigma-\tau} \ \mathrm{d}\sigma\Big)(\tau-A)^{-1} \ \mathrm{d}\tau\\
		&=\frac{1}{2\pi i}\int_{\partial\Omega}\varphi(\sigma)\Big(\frac{1}{2\pi i}\int_{\partial D}\frac{1}{\sigma-\tau}(\tau-A)^{-1} \ \mathrm{d}\tau\Big) \ \mathrm{d}\sigma\\
		&=\frac{1}{2\pi i}\int_{\partial\Omega}\varphi(\sigma)(\sigma-A)^{-1} \ \mathrm{d}\sigma
	\end{align*}
        and, likewise,
        \begin{align*}
            \gamma (\Phi_{\Omega}^{\ins}(\varphi^{*}))^{*}=\Big(\frac{1}{2\pi i}\int_{\partial\Omega}\varphi(\sigma)^{*}(\sigma-A)^{-1} \ \mathrm{d}\sigma\Big){}^{*}.
        \end{align*}
        So by definition of the harmonic functional calculus we infer that
	\begin{align*}
	\hat{\gamma}(K_{\Omega}^{\ins}(\varphi))=\gamma (\Phi_{\Omega}^{\ins}(\varphi))+\gamma (\Phi_{\Omega}^{\ins}(\varphi^{*}))^{*}=\int_{\partial\Omega}\varphi(\sigma)P_{\Omega}(\sigma,A) \ |\mathrm{d}\sigma|
	\end{align*}
	as desired.
\end{proof}
In particular, using Proposition \ref{integralRepresentationOperatorValuedDP} one readily verifies that
\begin{align*}
	\int_{\partial\Omega}P_{\Omega}(\sigma,A) \ |\mathrm{d}\sigma|=2.
\end{align*}

\subsection{Inclusion of $W(A)$ in $\Omega^{-}$}
If $\Omega$ is convex, then the positivity of the operator-valued double-layer potential encodes the inclusion of $W(A)$ in the closure $\Omega^{-}$. The next result, which we quote from \cite[Lemma 2.1]{badea2006convex}, can also be found between the lines in \cite[Section 1]{delyon1999generalization} and \cite[Section 3]{putinar2005skew}. 
\begin{proposition}\label{characterizationContainmentNumericalRange1}
	Assume that $\Omega$ is convex. The inclusion $W(A)\subseteq\Omega^{-}$ holds if and only if $P_{\Omega}(\sigma,A)\geq0$ for all $\sigma\in\partial\Omega$.
\end{proposition}
\begin{proof}
	For each point $\sigma$ on the boundary $\partial\Omega$ let $\Pi_{\sigma}$ be the unique closed half-plane in $\mathbb{C}$ with $n_{\Omega}(\sigma)$ as outward normal vector at $\sigma$. Convexity of $\Omega$ gives
	\begin{align*}
		\Omega^{-}=\bigcap_{\sigma\in\partial\Omega}\Pi_{\sigma}.
	\end{align*}
    For any $x,y\in H$ that satisfy $x=(\sigma-A)^{-1}y$ and $\|x\|=1$ the equality
    \begin{align*}
        \Re((\sigma-\langle Ax,x\rangle)^{*}n_{\Omega}(\sigma))=\pi\langle P(\sigma,A)y,y\rangle
    \end{align*}
    holds, implying that $\langle Ax,x\rangle\in\Pi_{\sigma}$ if and only if $\langle P(\sigma,A)y,y\rangle\geq0$. The result follows.
\end{proof}

In Proposition \ref{characterizationConvexity2} we have seen that convexity of $\Omega$ is encoded by the operator norm of $K_{\Omega}$. In a similar fashion, provided that $\Omega$ is convex, the inclusion of $W(A)$ in $\Omega^{-}$ is encoded by the operator norm of $\hat{\gamma}\circ K_{\Omega}^{\ins}$. The following result is a consequence of the integral representation of $\hat{\gamma}\circ K_{\Omega}^{\ins}$ presented in Proposition \ref{integralRepresentationOperatorValuedDP}.

\begin{proposition}\label{characterizationContainmentNumericalRange2}
	Assume that $\Omega$ is convex. The following assertions are equivalent:
		\begin{itemize}
			\item[{\normalfont(i)}] $\Omega^{-}$ contains $W(A)$,
			\item[{\normalfont(ii)}] $\hat{\gamma}\circ K_{\Omega}^{\ins}$ is a positive operator,
			\item[{\normalfont(iii)}] the operator norm of $\hat{\gamma}\circ K_{\Omega}^{\ins}$ equals $2$.
		\end{itemize}
\end{proposition}
\begin{proof}
 It follows from Proposition \ref{characterizationContainmentNumericalRange1} that $\Omega^{-}$ contains $W(A)$ if and only if $P_{\Omega}(\sigma,A)$ is positive for every point $\sigma$ on $\partial\Omega$. By Proposition \ref{integralRepresentationOperatorValuedDP} the latter is true if and only if $\tfrac{1}{2}\hat{\gamma}\circ K_{\Omega}^{\ins}$ is positive. Since $\tfrac{1}{2}\hat{\gamma}\circ K_{\Omega}^{\ins}$ is a unital operator between C*-algebras, it is well-known that $\tfrac{1}{2}\hat{\gamma}\circ K_{\Omega}^{\ins}$ is positive if and only if $\tfrac{1}{2}\hat{\gamma}\circ K_{\Omega}^{\ins}$ is contractive, see e.g.\ \cite[Chapter 2]{paulsen2002completely}.
\end{proof}
\begin{remark}
Assume that $\Omega$ is convex. If $\Omega^{-}$ contains $W(A)$, then the operator norm of $\hat{\gamma}\circ K_{\Omega}^{\ins}|\mathcal{A}(\Omega^{-})$ equals $2$ as well. It is an interesting question whether the converse of this statement is true.
\end{remark}

\section{Putinar--Sandberg's theorem}\label{sectionPutinar}
Let $A$ be a bounded operator on a Hilbert space $H$. Let $\Omega$ be a smoothly bounded open subset of $\mathbb{C}$ that contains the spectrum of $A$. Assume that $\Omega$ is connected. By Proposition \ref{decomposition} the mappings $\hat{\gamma}\circ K_{\Omega}^{\ins}$ and $\gamma$ agree on $\ker(K_{\Omega}|\mathcal{A}(\Omega^{-}))$. Putinar--Sandberg used this fact to prove the following result, see \cite[Theorem 3]{putinar2005skew}. The proof we provide here takes inspiration from the reasoning of Putinar--Sandberg, but circumvents the dilation theoretic arguments.

\begin{theorem}\label{strongBergerStampfli}
        Assume that $\Omega$ is convex. If the inclusion $W(A)\subseteq\Omega^{-}$ holds, then for any $f\in\mathcal{A}(\Omega^{-})$ with $\|f\|_{\infty}\leq1$ and $K_{\Omega}(f)=0$ one has $W(\gamma(f))\subseteq\mathbb{D}^{-}$.
\end{theorem}
\begin{proof}
		Take $f\in\mathcal{A}(\Omega^{-})$ with $\|f\|_{\infty}\leq1$ and $K_{\Omega}(f)=0$. Note that the spectrum of $\gamma(f)$ is contained in $\mathbb{D}$ by the maximum modulus principle. So by Proposition \ref{characterizationContainmentNumericalRange1} it suffices to prove that $P_{\mathbb{D}}(\tau,\gamma(f))\geq0$ for all $\tau\in\partial\mathbb{D}$.
        Let $\tau\in\partial\mathbb{D}$ and $r\in\interval[open right]{0}{1}$ be arbitrary. Recall from Proposition \ref{kernelK} that $\ker(K_{\Omega}|\mathcal{A}(\Omega^{-}))$ is a subalgebra of $\mathcal{A}(\Omega^{-})$. So for $k=1,2,\ldots$ we infer that $K_{\Omega}(f^{k})=0$ and therefore $\gamma(f^{k})=(\hat{\gamma}\circ K_{\Omega}^{\ins})(f^{k})$. Hence
	\begin{align*}
		P_{\mathbb{D}}(\tau,\gamma(rf))&=\frac{1}{\pi}\Re(\tau(\tau-\gamma(rf))^{-1})\\
		&=\frac{1}{\pi}\Re\Big(1+\sum_{k=1}^{\infty}\frac{r^{k}}{\tau^{k}}\gamma(f^{k})\Big)\\
		&=\frac{1}{\pi}\Re\Big(1+\sum_{k=1}^{\infty}\frac{r^{k}}{\tau^{k}}(\hat{\gamma}\circ K_{\Omega}^{\ins})(f^{k})\Big)\\
		&=\frac{1}{2\pi}\Re\Big((\hat{\gamma}\circ K_{\Omega}^{\ins})\Big(1+2\sum_{k=1}^{\infty}\frac{r^{k}}{\tau^{k}}f^{k}\Big)\Big)\\
		&=\frac{1}{2\pi}\Re((\hat{\gamma}\circ K_{\Omega}^{\ins})((\tau-rf)^{-1}(\tau+rf))).
	\end{align*}
	Since $P_{\Omega}(\sigma,A)\geq0$ for all $\sigma\in\partial\Omega$ by Proposition \ref{integralRepresentationOperatorValuedDP}, we obtain
	\begin{align*}
		P_{\mathbb{D}}(\tau,\gamma(rf))&=\frac{1}{2\pi}\int_{\partial\Omega}\Re\Big(\frac{\tau+rf(\sigma)}{\tau-rf(\sigma)}\Big)P_{\Omega}(\sigma,A) \ |\mathrm{d}\sigma|\geq0
	\end{align*}
	by Proposition \ref{characterizationContainmentNumericalRange1}. Letting $r\nearrow1$ completes the proof.
	\end{proof}
 Suppose that $f\in\mathcal{A}(\mathbb{D}^{-})$ and $\|f\|_{\infty}\leq1$. Recall from Proposition \ref{valuesDiskCaseK} that $K_{\mathbb{D}}(f)=f(0)$. In the case $\Omega=\mathbb{D}$ the result from Theorem \ref{strongBergerStampfli} was already discovered independently by Kato in 1965 and Berger--Stampfli in 1967 and states that, if $W(A)$ is contained in $\mathbb{D}^{-}$, then the same is true for $W(\gamma(f))$ provided that $f(0)=0$, see \cite[Theorem 5]{kato1965some} and \cite[Theorem 4]{berger1967mapping}. The proof of Kato was based on the Riesz--Herglotz representation theorem and the proof of Berger--Stampfli was based on dilation theory. In 2008 Drury analyzed the situation where $f(0)\neq0$ and proved that, if $W(A)$ is contained in $\mathbb{D}^{-}$, then $W(\gamma(f))$ is contained in the `teardrop region' described by the convex hull of $\mathbb{D}^{-}$ and the closed disk with center $f(0)$ and radius $1-|f(0)|^{2}$, see \cite[Theorem 2]{drury2008symbolic}. Alternative proofs of these results can be found in \cite{klaja2016mapping}.
 
\begin{corollary}\label{weakBergerStampfli}
 Assume that $\Omega$ is convex. If the inclusion $W(A)\subseteq\Omega^{-}$ holds, then for any $f\in\mathcal{A}(\Omega^{-})$ with $\|f\|_{\infty}\leq1$ and $K_{\Omega}(f)=0$ one has $\|\gamma(f)\|\leq2$.
\end{corollary}
\begin{proof}
	Take $f\in\mathcal{A}(\Omega^{-})$ with $\|f\|_{\infty}\leq1$ and $K_{\Omega}(f)=0$. Since $\|\gamma(f)\|\leq2w(\gamma(f))$, the result follows immediately from Theorem \ref{strongBergerStampfli}.
\end{proof}
In the case $\Omega=\mathbb{D}$ it follows from a result published in 1975 by Okubo--Ando that Corollary \ref{weakBergerStampfli} still holds without the assumption $f(0)=0$, see \cite[Theorem 2]{okubo1975constants}. Their proof is based on dilation theory.

\section{The role of $\hat{\gamma}\circ K_{\Omega}^{\mathrm{i}}$ in the literature on spectral constants}\label{sectionRole}
Let $A$ be a bounded operator on a Hilbert space $H$. Let $\Omega$ be a smoothly bounded operator on a Hilbert space $H$ that contains the spectrum of $A$. Assume that $\Omega$ is convex and that $\Omega^{-}$ contains $W(A)$. In this section we give an overview of the various appearances of $\hat{\gamma}\circ K_{\Omega}^{\ins}$ in the search for spectral constants of $A$ on $\Omega$. The strategies in the literature below boil down to finding a relation between the (unknown) norm of $\gamma$ and the (known) norm of $\hat{\gamma}\circ K_{\Omega}^{\ins}$.
\subsection{Delyon--Delyon (1999)}
Delyon--Delyon \cite{delyon1999generalization} proved that $K_{\Omega}^{\ins}$ is invertible as a self-mapping on $\mathcal{C}(\partial\Omega)$, which quickly leads to the estimate
\begin{align*}
\|\gamma(f)\|=\|\hat{\gamma}(f)\|=\|(\hat{\gamma}\circ K_{\Omega}^{\ins})((K_{\Omega}^{\ins})^{-1}(f))\|\leq 2\|(K_{\Omega}^{\ins})^{-1}(f)\|_{\infty}
\end{align*}
for all $f\in\mathcal{A}(\Omega^{-})$ with $\|f\|_{\infty}\leq1$. This estimate is sharp in the sense that it becomes an equality for constant functions. After taking suprema we obtain $\|\gamma\|\leq2\|(K_{\Omega}^{\ins})^{-1}|\mathcal{A}(\Omega^{-})\|$. In the same paper Delyon--Delyon estimated
\begin{align*}
    \|(K_{\Omega}^{\ins})^{-1}\|\leq\frac{1}{2}\Big(3+\Big(\frac{2\pi d_{\Omega}^{2}}{a_{\Omega}}\Big){}^{3}\Big),
\end{align*}
where $a_{\Omega}$ and $d_{\Omega}$ denote the area and diameter of the domain $\Omega$, respectively. 

\subsection{Crouzeix (2007)}
Crouzeix \cite{crouzeix2007numerical} used the decomposition from Proposition \ref{decomposition} and applied the triangle-inequality to obtain
\begin{align*}
    \|\gamma(f)\|=\|(\hat{\gamma}\circ K_{\Omega}^{\ins})(f)-\gamma(K_{\Omega}(f)^{*})^{*}\|\leq2+\|\gamma(K_{\Omega}(f)^{*})\|
\end{align*}
for all $f\in\mathcal{A}(\Omega^{-})$ with $\|f\|_{\infty}\leq1$. This estimate provides an alternative argument for Corollary \ref{weakBergerStampfli}. Remarkably, in the same paper Crouzeix estimated the operator norm of $\gamma(K_{\Omega}(f)^{*})$ by the absolute constant $9.08$.
\subsection{Crouzeix--Palencia (2017)}
Under the assumption that $\|\gamma\|>1$, Crouzeix--Palencia \cite{crouzeix2017numerical} considered for any $f\in\mathcal{A}(\Omega^{-})$ with $\|f\|_{\infty}\leq1$ the function
\begin{align*}
    h_{f}:=(\|\gamma\|^{2}+K_{\Omega}(f)^{*}f)^{-1}f
\end{align*}
and argued that the operator $1-\gamma(h_{f})^{*}(\hat{\gamma}\circ K_{\Omega}^{\ins})(f)$ is non-invertible. This readily implies that
\begin{align*}
    1\leq \|\gamma(h_{f})^{*}(\hat{\gamma}\circ K_{\Omega}^{\ins})(f)\|\leq 2\|\gamma\|\|h_{f}\|_{\infty}\leq\frac{2\|\gamma\|}{\|\gamma\|^{2}-1}
\end{align*}
and therefore $\|\gamma\|\leq1+\sqrt{2}$.
\subsection{Ransford--Schwenninger (2018)}
The authors of \cite{ransford2018remarks} and later \cite{ostermann2020abstract} presented an abstract framework and alternative proof for Crouzeix--Palencia's result. Using the C*-identity and triangle-inequality they deduced for any $f\in\mathcal{A}(\Omega^{-})$ with $\|f\|_{\infty}\leq1$ that
\begin{align*}
\|\gamma(f)\|^{4}&=\|\gamma(f)^{*}\gamma(f)\gamma(f)^{*}\gamma(f)\|\\
&=\|\gamma(f)^{*}(\hat{\gamma}\circ K_{\Omega}^{\ins})(f)\gamma(f)^{*}\gamma(f)-\gamma(fK_{\Omega}(f)^{*}f)^{*}\gamma(f)\|\\
&\leq2\|\gamma(f)\|^{3}+\|\gamma(fK_{\Omega}(f)^{*}f)\|\|\gamma(f)\|.
\end{align*}
This establishes $\|\gamma\|^{4}\leq 2\|\gamma\|^{3}+\|\gamma\|^{2}$ and therefore $\|\gamma\|\leq 1+\sqrt{2}$. 

\section{A spectral constant}\label{sectionSpectral}
Let $A$ be a bounded operator on a Hilbert space $H$. Let $\Omega$ be a smoothly bounded open subset of $\mathbb{C}$ that contains the spectrum of $A$. Assume that $\Omega$ is connected. An element $f_{0}$ of $\mathcal{A}(\Omega^{-})$ with unit length is called an \textit{extremal function} for $A$ on $\Omega$ if $\gamma$ attains its operator norm at $f_{0}$. Similarly, an element $x_{0}$ of $H$ with unit length is called an \textit{extremal vector} for $A$ on $\Omega$ associated to $f_{0}$ if $\gamma(f_{0})$ attains its operator norm at $x_{0}$. If $f_{0}$ and $x_{0}$ form such an \textit{extremal pair}, then there exists a probability Radon measure $\mu_{0}$ on $\partial\Omega$, called an \textit{extremal measure} for $A$ on $\Omega$ associated to $f_{0}$ and $x_{0}$, that satisfies
\begin{align*}
\int_{\partial\Omega}\varphi \ \mathrm{d}\mu_{0}=\langle \hat{\gamma}(\varphi)x_{0},x_{0}\rangle
\end{align*}
for all $\varphi\in\mathcal{H}(\Omega^{-})$, see \cite[Theorem 2.1]{bickel2020crouzeix}. Extremal measures obey the formula
\begin{align*}
    (\|\gamma\|-1)\int_{\partial\Omega}f_{0} \ \mathrm{d}\mu_{0}=0,
\end{align*}
see \cite[Theorem 5.1]{caldwell2018some}. More abstract versions of these results can be found in \cite{schwenninger2024abstract}.

Assume that there exists an extremal pair $f_{0}$ and $x_{0}$ for $A$ on $\Omega$ and consider the constant
\begin{align*}
\rho(f_{0},x_{0}):=\langle \hat{\gamma}(\Re(K_{\Omega}(f_{0})^{*}f_{0}))x_{0},x_{0}\rangle.
\end{align*}
Using an extremal measure $\mu_{0}$ we see that
\begin{align*}
    \rho(f_{0},x_{0})=\int_{\partial\Omega}\Re(K(f_{0})^{*}f_{0}) \ \mathrm{d}\mu_{0}
\end{align*}
and therefore
\begin{align*}
|\rho(f_{0},x_{0})|\leq\int_{\partial\Omega}|\Re(K_{\Omega}(f_{0})^{*}f_{0})| \ \mathrm{d}\mu_{0}\leq\mu_{0}(\partial\Omega)\|K_{\Omega}(f_{0})\|_{\infty}\|f_{0}\|_{\infty}\leq\|K_{\Omega}\|.
\end{align*}
In particular, if $\Omega$ is convex, then $|\rho(f_{0},x_{0})|\leq1$ by Proposition \ref{characterizationConvexity2}. Adopting techniques from \cite{schwenninger2024abstract} we now prove the following result.

\begin{theorem}\label{spectralConstant}
Assume that $\Omega$ is convex. Suppose that there exists an extremal pair $f_{0}$ and $x_{0}$ for $A$ on $\Omega$. If $W(A)$ is contained in $\Omega^{-}$, then $1+\sqrt{1-\rho(f_{0},x_{0})}$ is a spectral constant for $A$ on $\Omega$.
\end{theorem}
\begin{proof}
By Proposition \ref{characterizationContainmentNumericalRange2} we have $\|\hat{\gamma}\circ K_{\Omega}^{\ins}\|=2$ and therefore
\begin{align*}
\rho(f_{0},x_{0})+\|\gamma\|^{2}&=\Re(\langle \gamma(f_{0})x_{0},\gamma(K_{\Omega}(f_{0})^{*})^{*}x_{0}\rangle)+\Re(\langle \gamma(f_{0})x_{0},\gamma(f_{0})x_{0}\rangle)\\
&=\Re(\langle\gamma(f_{0})x_{0},(\hat{\gamma}\circ K_{\Omega}^{\ins})(f_{0})x_{0}\rangle)\\
&\leq\|(\hat{\gamma}\circ K_{\Omega}^{\ins})(f_{0})x_{0}\|\|\gamma\|\leq2\|\gamma\|.
\end{align*}
Since $\rho(f_{0},x_{0})\leq1$, it follows that $\|\gamma\|\leq 1+\sqrt{1-\rho(f_{0},x_{0})}$ as desired.
\end{proof}
Assume that $H$ is finite-dimensional and that $\Omega$ is convex. In this case there always exists an extremal pair $f_{0}$ and $x_{0}$ for $A$, see \cite[Theorem 2.1]{crouzeix2004bounds}. Also assume that $\Omega^{-}$ contains $W(A)$. Since $\rho(f_{0},x_{0})\geq-1$, we find that $\|\gamma\|\leq1+\sqrt{2}$ by Theorem \ref{spectralConstant}. By a Krylov subspace argument the same bound holds in the infinite-dimensional case, see \cite[Theorem 2]{crouzeix2007numerical}. In particular, we recover Crouzeix--Palencia's result. Interestingly, if it is true that $\rho(f_{0},x_{0})\geq0$, then Crouzeix's conjecture follows.

Finally, suppose that $\Omega$ is an open disk. In this case $K_{\Omega}(f_{0})$ is constant by Proposition \ref{valuesDiskCaseK}. It follows that $(\|\gamma\|-1)\rho(f_{0},x_{0})=0$ and therefore $\|\gamma\|\leq2$ by Theorem \ref{spectralConstant}.

\begin{remark}
Let $C$ be any convex subset of $\mathbb{C}$ with non-empty interior. In the interesting article \cite{malman2024double} it was shown that the so-called analytic configuration $a(C)$, which in our notation with $C=\Omega^{-}$ is given by
\begin{align*}
    a(\Omega^{-})=\inf_{\lambda\in\mathbb{C}}\sup_{\substack{f\in\mathcal{A}(\Omega^{-}) \\ \|f\|_{\infty}\leq1}}\|K_{\Omega}(f)-\lambda\|_{\infty},
\end{align*}
is strictly smaller than $1$. Since for any $\lambda\in\mathbb{C}$ the estimate
\begin{align*}
    |\rho(f_{0},x_{0})|=\Big|\Re\Big(\int_{\partial\Omega}(K_{\Omega}(f_{0})-\lambda)^{*}f_{0} \ \mathrm{d}\mu_{0}\Big)\Big|\leq\|K_{\Omega}(f_{0})-\lambda\|_{\infty}
\end{align*}
holds (under the natural assumption $\|\gamma\|\neq1$), this readily implies that
\begin{align*}
    |\rho(f_{0},x_{0})|\leq a(\Omega^{-})<1.
\end{align*}
Note that in \cite{malman2024double} the authors also provide the spectral constant $1+\sqrt{1+a(C)}$ for $A$ on any $C$ that contains $W(A)$. While for $\Omega^{-}$ this spectral constant seems slightly larger than ours from Theorem \ref{spectralConstant}, they can consider more general convex domains $C$.
\end{remark}

\appendix
\section{Approximation by smooth sets}
We say that an open subset $\Omega$ of $\mathbb{C}$ is \textit{smoothly bounded} if $\Omega$ is bounded and the boundary $\partial\Omega$ is an embedded smooth submanifold of $\mathbb{C}$.

Given a non-empty subset $X$ of $\mathbb{C}$ and an element $z$ in $\mathbb{C}$, we write $\dis(X,z)$ for the infimum of all Euclidian distances between points in $X$ and the point $z$. The mapping
\begin{align*}
    z\mapsto\dis(X,z),\qquad z\in\mathbb{C}
\end{align*}
is called the \textit{distance function} of $X$. Any distance function is Lipschitz continuous with Lipschitz constant $1$. If $X$ is convex, then the distance function of $X$ is convex as well.

\begin{lemma}\label{smoothNeighbourhoods}
Let $X$ be a bounded non-empty subset of $\mathbb{C}$ and $E$ an open neighbourhood of the closure $X^{-}$. There exist smoothly bounded open subsets $\Omega_{1},\Omega_{2},\ldots$ of $\mathbb{C}$ containing $X^{-}$ such that $\Omega_{n+1}^{-}\subseteq\Omega_{n}$ and $\Omega_{n}^{-}\subseteq E$ for all $n\in\mathbb{N}$ and, furthermore,
\begin{align*}
X^{-}=\bigcap_{n\in\mathbb{N}}\Omega_{n}^{-}.
\end{align*}
If $X$ is convex, then $\Omega_{1},\Omega_{2},\ldots$ may be taken convex as well.
\end{lemma}
\begin{proof}
For each $n\in\mathbb{N}$ we define the neighbourhood
\begin{align*}
    N_{n}:=\Big\{z\in\mathbb{C}:\dis(X,z)<\frac{\varepsilon}{n}\Big\},
\end{align*}
where $\varepsilon>0$ is chosen such that $N_{1}\subseteq E$. Fix a number
\begin{align*}
    0<s_{n}<\frac{1}{2}\Big(\frac{\varepsilon}{n}-\frac{\varepsilon}{n+1}\Big)
\end{align*}
and define $\theta_{n}\colon\mathbb{C}\to\mathbb{R}$ by
\begin{align*}
\theta_{n}(z):=\frac{1}{s_{n}^{2}}\theta\Big(\frac{z}{s_{n}}\Big)
\end{align*}
for $z\in\mathbb{C}$, where $\theta\colon\mathbb{C}\to\mathbb{R}$ is the standard mollifier on the plane. Now consider the smooth positive function $\psi_{n}\colon\mathbb{C}\to\mathbb{R}$ given by
\begin{align*}
    \psi_{n}(z):=\int_{\mathbb{C}}\theta_{n}(z-\lambda)\dis(N_{n+1},\lambda) \ \mathrm{d}\lambda
\end{align*}
for $z\in\mathbb{C}$. By Sard's theorem, see e.g.\ \cite[Theorem 6.10]{lee2012introduction}, this function has a regular value
\begin{align*}
s_{n}<t_{n}<\frac{1}{2}\Big(\frac{\varepsilon}{n}-\frac{\varepsilon}{n+1}\Big)
\end{align*}
so that the boundary of
\begin{align*}
    \Omega_{n}:=\{z\in\mathbb{C}:\psi_{n}(z)<t_{n}\}
\end{align*}
is an embedded smooth submanifold. 

Let us first prove that $N_{n+1}^{-}\subseteq\Omega_{n}$. If $z\in N_{n+1}^{-}$, then $\dis(N_{n+1},z)=0$ and therefore
\begin{align*}
\psi_{n}(z)&=\int_{\mathbb{C}}\theta_{n}(z-\lambda)|\dis(N_{n+1},z)-\dis(N_{n+1},\lambda)| \ \mathrm{d}\lambda\\
&\leq\int_{\mathbb{C}}\theta_{n}(z-\lambda)|z-\lambda| \ \mathrm{d}\lambda\leq s_{n}<t_{n},
\end{align*}
which gives $z\in \Omega_{n}$.

Next we show that $\Omega_{n}^{-}\subseteq N_{n}$. If $z\in\Omega_{n}^{-}$, then $\psi_{n}(z)\leq t_{n}$ and therefore
\begin{align*}
\dis(X,z)&\leq\frac{\varepsilon}{n+1}+\dis(N_{n+1},z)-\psi_{n}(z)+t_{n}\\
&\leq\frac{\varepsilon}{n+1}+\int_{\mathbb{C}}\theta_{n}(z-\lambda)|\dis(N_{n+1},z)-\dis(N_{n+1},\lambda)| \ \mathrm{d}\lambda+t_{n}\\
&\leq\frac{\varepsilon}{n+1}+s_{n}+t_{n}<\frac{\varepsilon}{n},
\end{align*}
which gives $z\in N_{n}$.

It is clear that $\Omega_{1},\Omega_{2},\ldots$ satisfy the desired properties. If $X$ is convex, then $N_{n+1}$ is convex, which in turn implies that the distance function of $N_{n+1}$ is convex. Hence $\psi_{n}$ is convex (being the mollification of a convex function) and $\Omega_{n}$ is convex (being the sublevel set of a convex function).
\end{proof}
For any open subset $E$ of $\mathbb{C}$ we write $\mathscr{C}_{E}$ for the space of all non-empty compact subsets of $E$ and equip it with the Hausdorff metric. It is well-known that, if a sequence $(C_{n})_{n\in\mathbb{N}}$ in $\mathscr{C}_{E}$ is nested in the sense that $C_{1}\supseteq C_{2}\supseteq\ldots$, then it converges to its intersection in $\mathscr{C}_{E}$, see e.g.\ \cite[Lemma 1.8.2]{schneider2013convex}.

\begin{proposition}\label{smoothApproximation}
Let $A$ be an operator on a Hilbert space $H$. Suppose that $X$ is bounded subset of $\mathbb{C}$ such that $X^{-}$ contains the spectrum of $A$. The following statements are equivalent:
\begin{itemize}
			\item[{\normalfont(i)}] $\kappa$ is a spectral constant for $A$ on $X$,
			\item[{\normalfont(ii)}] $\kappa$ is a spectral constant for $A$ on every smoothly bounded open neighbourhood $\Omega$ of the spectrum of $A$ for which $\Omega^{-}$ contains $X$.
		\end{itemize}
If $X$ is convex, then the previous statements are also equivalent to the following statement:
\begin{itemize}
    \item[{\normalfont(iii)}] $\kappa$ is a spectral constant for $A$ on every smoothly bounded open neighbourhood $\Omega$ of the spectrum of $A$ for which $\Omega$ is convex and $\Omega^{-}$ contains $X$.
\end{itemize}
\end{proposition}
\begin{proof}
    It is clear that (i) implies (ii). To prove the converse, consider any rational function $f\colon X^{-}\to\mathbb{C}$ with poles off $X^{-}$. Let $E$ be the complement of the finitely many poles of $f$ in $\mathbb{C}$. By Lemma \ref{smoothNeighbourhoods} there exist smoothly bounded open subsets $\Omega_{1},\Omega_{2},\ldots$ of $\mathbb{C}$ containing $X^{-}$ such that $\Omega_{n+1}^{-}\subseteq\Omega_{n}$ and $\Omega_{n}^{-}\subseteq E$ for all $n\in\mathbb{N}$ and, furthermore,
    \begin{align*}
        X^{-}=\bigcap_{n\in\mathbb{N}}\Omega_{n}^{-}.
    \end{align*}
    We infer that $(\Omega_{n}^{-})_{n\in\mathbb{N}}$ converges to $X^{-}$ in $\mathscr{C}_{E}$. For each positive integer $n$ let $f_{n}$ be the unique rational extension of $f$ to $\Omega_{n}^{-}$. It follows that
    \begin{align*}
        \|f(A)\|=\lim_{n\to\infty}\|f_{n}(A)\|\leq\lim_{n\to\infty}\kappa\cdot\sup_{z\in\Omega_{n}}|f_{n}(z)|=\kappa\cdot\sup_{z\in X}|f(z)|.
    \end{align*}
    If $X^{-}$ is convex, then (i) and (iii) are also equivalent because $\Omega_{1},\Omega_{2},\ldots$ may be taken convex in this case.
\end{proof}

\subsection*{Data availability}
Not applicable.

\subsection*{Declaration of interest}
The authors declare that they have no financial interest or conflict to disclose. 

\bibliographystyle{alpha}
\bibliography{Bibliography}

\end{document}